\documentclass[12pt]{article}
\usepackage{latexsym,amsfonts,amssymb}
\usepackage[dvips]{color}
\usepackage{colordvi,multicol}
\usepackage[marginal]{footmisc}

\def \cal{\mathcal}

\newtheorem{thm}{Theorem}[section]

\newtheorem{pro}[thm]{Proposition}
\newtheorem{defi}[thm]{Definition}
\newtheorem{rem}[thm]{Remark}
\newtheorem{exa}[thm]{Example}

\date{}

\begin{document}
\title{\bf  Weak Random Periodic Solutions of Random Dynamical
Systems}
\author{}

\maketitle

\centerline{Wei Sun} \centerline{\small Department of Mathematics
and Statistics}
\centerline{\small Concordia University} \centerline{\small
Montreal, H3G 1M8, Canada} \centerline{\small E-mail:
wei.sun@concordia.ca}
\vskip 1.5cm
\centerline{Zuo-Huan Zheng} \centerline{\small Institute of Applied Mathematics} \centerline{\small
Academy of Mathematics and Systems Science}\centerline{\small
Chinese Academy of Sciences}\centerline{\small
Beijing, 100190, China}\centerline{\small
and}
\centerline{\small
School of Mathematical Sciences}
\centerline{\small University of Chinese Academy of Sciences} \centerline{\small
Beijing, 100049, China}
\centerline{\small E-mail:
zhzheng@amt.ac.cn}
\vskip 0.2cm


\vskip 0.5cm \noindent{\bf Abstract:}\quad  We
first introduce the concept of weak random periodic solutions
of random dynamical systems. Then, we discuss the existence of
such periodic solutions.  Further, we introduce the definition of
weak random periodic measures and study their relationship with
weak random periodic solutions. In particular, we establish the
existence of invariant measures of random dynamical systems by
virtue of their weak random periodic solutions. Finally, we use concrete examples to
illustrate the weak random
periodic phenomena of dynamical systems induced by random and stochastic differential equations.

\vskip 0.5cm


\vskip 0.3cm









\section{Introduction}
Periodic solutions are a very active research topic of the
qualitative theory of ordinary differential equations. Given a
dynamical system, it is important to investigate the existence,
number and positions of periodic solutions as well as the behavior
of their nearby trajectories. For example, for polynomial vector fields in the plane, it is essential
to derive an upper bound of the number of limit cycles
and discuss their
relative positions. This is the second part of Hilbert's 16th
problem.  Another example is the famous Poincar\'e-Bendixson
theorem. It plays a fundamental role in the qualitative theory of
differential equations in the plane because it provides a useful
method to check the existence of periodic solutions and
to find their positions (cf. \cite{P,B, NS, H}). In the past
fifty years, a lot of progress has been made for periodic
solutions and the global structure of dynamical systems (cf.
\cite{Y1,Y,F, WY,ZXZ,LW, MWSY, LLW}).

When we study random dynamical systems, it is natural to consider the counterparts of fixed points and periodic solutions. In the literature, the counterparts are called stationary solutions and
random periodic solutions, respectively. Stationary solutions have
attracted lots of attention and a series of results have been
obtained (cf. \cite{A,ZZ, MZZ,LZ}).  In \cite{ZZ2}, Zhao and Zheng
introduced for the first time the concept of random periodic
solutions and gave a sufficient condition for their existence. In
\cite{FZ}, Feng and Zhao introduced random periodic measures and
discussed the close relationship between random periodic solutions
and random periodic measures. We call the
reader's attention to \cite{FZZ,FZ6,FWZ3,FLZ1,D1,D2} for other recent works
on random periodic solutions.

Note that the period $T$ of the random periodic solutions
introduced in \cite{ZZ2} is deterministic and uniform for all random
paths $\omega$. However, for many random
dynamical systems induced by random or stochastic differential equations, the solutions exhibit some periodic behaviors while
the periods depend on $\omega$. To deal with these phenomena, we
 introduce in this paper the novel concept of weak random
periodic solutions. It is easy to see that any random periodic solution is a weak random periodic solution. But, in general, a weak random
periodic solution might not be a random periodic solution.

The remainder of this paper is organized as follows. In Section 2, we introduce the concept of weak random periodic solutions and present a useful criterion for their existence. In Section 3, we give the definition of
weak random periodic measures and show that the existence of weak
random periodic solutions implies the existence of weak random
periodic measures. Further, we establish the existence of
invariant measures for random dynamical systems by virtue of their
weak random periodic solutions.  In  Section 4, we use concrete examples  to illustrate the weak random
periodic phenomena of dynamical systems induced by random and stochastic differential equations.

\section{Definition and existence of weak random periodic solutions}

First, let us recall the concepts of fixed point and periodic
solution. Let $E$ be a Polish space with  Borel
$\sigma$-algebra ${\cal B}(E)$. For a
deterministic dynamical system $\Psi:\mathbb{R}\times E\rightarrow E$, a fixed point is a point
$x\in E$ such that
$$
\Psi(t)x=x,\ \ \ \ \forall t\in\mathbb{R}.
$$
A periodic solution with period $T>0$ is a ${\cal B}(\mathbb{R})$-measurable function $Y: \mathbb{R}\rightarrow
E$ such that
$$
\Psi(t)Y(s)=Y(t+s),\ \ Y(t+T)=Y(t),\ \ \ \ \forall
s,t\in\mathbb{R}.
$$

Suppose $\Psi:\mathbb{R}\times
\Omega\times E\rightarrow E$ is a measurable random dynamical
system on $(E,{\cal B}(E))$  over a metric dynamical system
$(\Omega,{\cal F},P,(\theta_t)_{t\in\mathbb{R}})$. Then,  for  $\omega\in\Omega$,
\begin{equation}\label{1}
\Psi(0,\omega)={\rm id}_E,\ \ \Psi(t+s,\omega)=\Psi(t,\theta_s\omega)\Psi(s,\omega),\ \ \ \ \forall s,t\in\mathbb{R}.
\end{equation}
A stationary solution (cf. \cite{A}) of $\Psi$ is a random variable $Y: {\Omega}\rightarrow E$ such
that for almost all $\omega\in\Omega$,
$$
\Psi(t,\omega)Y(\omega)=Y(\theta_t\omega),\ \ \ \ \forall
t\in\mathbb{R}.
$$
A random periodic solution with period $T>0$ (see \cite{ZZ2} and \cite{FZ}) is a ${\cal B}(\mathbb{R})\times{\cal
F}$-measurable function $Y:\mathbb{R}\times\Omega
\rightarrow E$ such that for almost all $\omega\in\Omega$,
$$
\Psi(t,\theta_s\omega)Y(s,\omega)=Y(t+s,\omega),\ \
Y(s+T,\omega)=Y(s,\theta_T\omega),\ \ \ \ \forall s,t\in\mathbb{R}.
$$

The period $T$ in the above definition is non-random. For many applications, this is not satisfactory. Here is a simple example.  Suppose $X(\varpi)$ is a positive random variable. Consider the following random differential equation (RDE):
\begin{equation}\label{KKLL}
\frac{d^2x(t)}{dt^2}=\sin(X(\varpi)(t+s)),
\end{equation}
where $s\in\mathbb{R}$. The periodic solution of (\ref{KKLL})  is given by
$$
x(t)=\frac{-\sin(X(\varpi)( t+s))}{X(\varpi)^2}+r,\ \ \ \ r\in\mathbb{R},
$$
whose period $T\varpi=\frac{2\pi}{X(\varpi)}$ is random.

To deal with the phenomenon of random periods, we now introduce the concept of weak random periodic solution of a random dynamical system.
\begin{defi}\label{defi2}
    A weak random periodic solution of $\Psi$ is a pair of measurable maps $Y: \mathbb{R}\times \Omega\rightarrow E$ and $T:\Omega\rightarrow(0,\infty)$ such that for almost all $\omega\in\Omega$,
    \begin{equation}\label{3}
    \Psi(t,\theta_s\omega)Y(s,\omega)=Y(t+s,\omega),\ \ Y(s+T\omega,\theta_{-T\omega}\omega)=Y(s,\omega),\ \ \ \ \forall s,t\in\mathbb{R}.
    \end{equation}
\end{defi}

Obviously, if $T$ is a constant map then the weak random periodic solution is reduced to the random periodic solution. For the existence of weak random periodic solutions, we have the
following useful criterion.
\begin{pro}\label{thm2}
If there exist measurable maps $Y_0:\Omega\rightarrow E$ and $T:\Omega\rightarrow(0,\infty)$ such that for almost all $\omega\in\Omega$,
\begin{equation}\label{26}
Y_0(\omega)=\Psi(T\omega,\theta_{-T\omega}\omega)Y_0(\theta_{-T\omega}\omega),
\end{equation}
Then, the random dynamical system $\Psi$ has a weak random periodic solution.
\end{pro}

\noindent {\bf Proof.}\ \ For $\omega\in\Omega$, define $Y(0,\omega)=Y_0(\omega)$ and
\begin{equation}\label{2}
Y(t,\omega):=\Psi(t,\omega)Y(0,\omega), \ \ t\in \mathbb{R}.
\end{equation}
Then, by (\ref{1}), (\ref{26}) and (\ref{2}), we obtain that for almost all $\omega\in\Omega$,
\begin{eqnarray*}
\Psi(t,\theta_s\omega)Y(s,\omega)&=&\Psi(t,\theta_s\omega)\Psi(s,\omega)Y(0,\omega)\\
&=&\Psi(t+s,\omega)Y(0,\omega)\\
&=&Y(t+s,\omega),\ \ \ \ \forall s,t\in\mathbb{R},
\end{eqnarray*}
and
\begin{eqnarray*}
Y(s,\omega)&=&\Psi(s,\omega)Y(0,\omega)\\
&=&\Psi(s,\omega)\Psi(T\omega,\theta_{-T\omega}\omega)Y(0,\theta_{-T\omega}\omega)\\
&=&\Psi(s+T\omega,\theta_{-T\omega}\omega)Y(0,\theta_{-T\omega}\omega)\\
&=&Y(s+T\omega,\theta_{-T\omega}\omega),\ \ \ \ \forall s\in \mathbb{R}.
\end{eqnarray*}
Therefore, $(Y,T)$ is a weak random periodic solution of $\Psi$.\hfill\fbox\\

We next consider the weak random periodic solution of a stochastic semi-flow. Denote
$$
\Delta=\{(t,s)\in \mathbb{R}^2:s\le t\}.
$$
Let $\varphi:\Delta\times\Omega\times E \rightarrow E$ be a stochastic semi-flow. Then, for  $\omega\in\Omega$,
\begin{equation}\label{1111}
\varphi(t,s,\omega)=\varphi(t,u,\omega)\circ\varphi(u,s,\omega),\ \ \ \ \forall s\le u\le t,
\end{equation}
and
$$
\varphi(s,s,\omega)={\rm id}_E,\ \ \ \ \forall s\in \mathbb{R}.
$$

\begin{defi}\label{defi22}
    A weak random periodic solution of $\varphi$ is a pair of measurable maps $Y: \mathbb{R}\times \Omega\rightarrow E$ and $T:\Omega\rightarrow(0,\infty)$ such that for almost all $\omega\in\Omega$,
    \begin{equation}\label{3333}
   \varphi(t,s,\omega)Y(s,\omega)=Y(t,\omega),\ \ Y(s+T\omega,\theta_{-T\omega}\omega)=Y(s,\omega),\ \ \ \ \forall s\le t.
    \end{equation}
\end{defi}

\section{Weak random periodic measures and invariant measures}\setcounter{equation}{0}

Let $\Psi$ be a measurable random dynamical
system. Define
$$
\Upsilon_t(\omega,x)=(\theta_t\omega,\Psi(t,\omega)x),\ \ \ \ \omega\in\Omega,\,x\in E,\, t\in\mathbb{R}.
$$
Denote by ${\cal P}(\Omega\times E)$ the set of all probability measures on $(\Omega\times E, \cal F\otimes\cal B(E))$.

\begin{defi}
    A weak random periodic probability measure of $\Psi$ is a pair of measurable maps $\mu: \mathbb{R}\times\Omega\rightarrow{\cal P}(\Omega\times E)$ and $T:\Omega\rightarrow(0,\infty)$ such that for almost all $\omega\in\Omega$,
    $$
    {\Upsilon}_t\mu(s,\omega)=\mu({t+s},\omega), \ \
    \mu(s+T\omega,\theta_{-T\omega}\omega)=\mu(s,\omega), \ \ \ \ \forall s,t\in\mathbb{R}.
    $$
\end{defi}

\begin{thm}\label{invar}
    If a random dynamical system $\Psi:\mathbb{R}\times\Omega\times E\rightarrow E$ has a weak random periodic solution $Y:\mathbb{R}\times\Omega\rightarrow E$
    and $T:\Omega\rightarrow(0,\infty)$, then it has a weak random periodic probability measure. Additionally, if for almost all $\omega\in\Omega$,
\begin{equation}\label{zsw}
T\omega=T(\theta_s\omega),\ \ \ \ \forall s\in \mathbb{R},
\end{equation}
then $\Psi$ has an invariant probability measure whose random factorization is supported by
$$
L^{\omega}:=\{Y(s,\theta_{-s}\omega):s\in[0,T\omega)\}.
$$
\end{thm}
\noindent {\bf Proof.}\ \ For $s\in\mathbb{R}$ and $\omega\in\Omega$, define
\begin{eqnarray*}
\mu(s,\omega)(A)=\delta_{Y(s,\omega)}(A_{\theta_s\omega}),\ \ \ \ A\in \cal F\otimes\cal B(E),
\end{eqnarray*}
where $A_{\omega}$ is the $\omega$-section of $A$. Then,  $\mu(s,\omega)\in {\cal P}(\Omega\times E)$.

We have
\begin{eqnarray*}
(\Upsilon^{-1}_t(A))_{\omega}=\{x:(\theta_t\omega,\Psi(t,\omega)x)\in A\}=\{x:\Psi(t,\omega)x\in A_{\theta_t\omega}\}=\Psi^{-1}(t,\omega)A_{\theta_t\omega}.
\end{eqnarray*}
Then, by (\ref{3}), we obtain that for almost all $\omega\in\Omega$,
\begin{eqnarray*}
{\Upsilon}_t\mu(s,\omega)(A)
&=&\mu(s,\omega)({\Upsilon}^{-1}_t(A))\nonumber\\
&=&\delta_{Y(s,\omega)}(({\Upsilon}^{-1}_t(A))_{\theta_s\omega})\nonumber\\
&=&\delta_{\Psi(t,\theta_s\omega)Y(s,\omega)}(A_{\theta_{t+s}\omega})\nonumber\\
&=&\delta_{Y(t+s,\omega)}(A_{\theta_{t+s}\omega})\nonumber\\
&=&\mu({t+s},\omega)(A),
\end{eqnarray*}
and
\begin{eqnarray*}
\mu({s+T\omega},\theta_{-T\omega}\omega)(A)
&=&\delta_{Y({s+T\omega},\theta_{-T\omega}\omega)}(A_{\theta_{s}\omega})\nonumber\\
&=&\delta_{Y(s,\omega)}(A_{\theta_{s}\omega})\nonumber\\
&=&\mu({s},\omega)(A).
\end{eqnarray*}
Thus, $\mu$ is a weak random periodic probability measure of
$\Psi$.

For $A\in \cal F\otimes\cal B(E)$, define
$$
\tilde{\mu}(A):=\int_{\Omega}\frac{1}{T\omega}\int_0^{T\omega}\mu(s,\omega)(A)dsP(d\omega).
$$
By (\ref{3}),  (\ref{zsw}) and the measure preserving property of $\{\theta_t\}$, we get
\begin{eqnarray*}
\tilde{\mu}(A)&=&\int_{\mathbb{R}}\int_{\Omega}\frac{\delta_s{([0,T\omega])}\cdot\delta_{Y(s,\omega)}(A_{\theta_{s}\omega})}{T\omega}P(d\omega)ds\\
&=&\int_{\mathbb{R}}\int_{\Omega}\frac{\delta_s{([0,T\omega])}\cdot\delta_{Y(s+T\omega,\theta_{-T\omega}\omega)}(A_{\theta_{s}\omega})}{T\omega}P(d\omega)ds\\
&=&\int_{\mathbb{R}}\int_{\Omega}\frac{\delta_s{([0,T\omega])}\cdot\delta_{Y(s+T\omega,\theta_{-s-T\omega}\omega)}(A_{\omega})}{T\omega}P(d\omega)ds\\
&=&\int_{\Omega}\int_{\mathbb{R}}\frac{\delta_s{([0,T\omega])}\cdot\delta_{Y(s+T\omega,\theta_{-s-T\omega}\omega)}(A_{\omega})}{T\omega}dsP(d\omega)\\
&=&\int_{\Omega}\int_{\mathbb{R}}\frac{\delta_s{([T\omega,2T\omega])}\cdot\delta_{Y(s,\theta_{-s}\omega)}(A_{\omega})}{T\omega}dsP(d\omega)\\
&=&\int_{\mathbb{R}}\int_{\Omega}\frac{\delta_s{([T\omega,2T\omega])}\cdot\delta_{Y(s,\theta_{-s}\omega)}(A_{\omega})}{T\omega}P(d\omega)ds\\
&=&\int_{\mathbb{R}}\int_{\Omega}\frac{\delta_s{([T\omega,2T\omega])}\cdot\delta_{Y(s,\omega)}(A_{\theta_{s}\omega})}{T\omega}P(d\omega)ds\\
&=&\int_{\Omega}\frac{1}{T\omega}\int_{T\omega}^{2T\omega}\mu(s,\omega)(A)dsP(d\omega).
\end{eqnarray*}
Repeating this argument, we can show that
\begin{eqnarray*}
\tilde{\mu}(A)=\int_{\Omega}\frac{1}{T\omega}\int_{k(T\omega)}^{(k+1)(T\omega)}\mu(s,\omega)(A)dsP(d\omega),\ \ \ \ \forall k\in\mathbb{N},
\end{eqnarray*}
which implies that
\begin{eqnarray}\label{4}
\tilde{\mu}(A)&=&\int_{\Omega}\lim_{N\rightarrow\infty}\frac{1}{N}\int_0^{N}\delta_{Y(s,\omega)}(A_{\theta_{s}\omega})dsP(d\omega)\nonumber\\
&=&\lim_{N\rightarrow\infty}\frac{1}{N}\int_0^{N}\int_{\Omega}\delta_{Y(s,\omega)}(A_{\theta_{s}\omega})P(d\omega)ds.
\end{eqnarray}

By (\ref{3}), (\ref{4})  and the measure preserving property of $\{\theta_t\}$, we obtain  that
\begin{eqnarray*}
\Upsilon_t\tilde{\mu}(A)&=&\tilde{\mu}(\Upsilon^{-1}_t(A))\\
&=&\lim_{N\rightarrow\infty}\frac{1}{N}\int_0^{N}\int_{\Omega}\delta_{Y(s,\omega)}((\Upsilon^{-1}_t(A))_{\theta_{s}\omega})P(d\omega)ds\\
&=&\lim_{N\rightarrow\infty}\frac{1}{N}\int_0^{N}\int_{\Omega}\delta_{\Psi(t,\theta_s\omega)Y(s,\omega)}(A_{\theta_{t+s}\omega})P(d\omega)ds\\
&=&\lim_{N\rightarrow\infty}\frac{1}{N}\int_0^{N}\int_{\Omega}\delta_{Y(t+s,\omega)}(A_{\theta_{t+s}\omega})P(d\omega)ds\\
&=&\lim_{N\rightarrow\infty}\frac{1}{N}\int_0^{N}\int_{\Omega}\delta_{Y(s,\omega)}(A_{\theta_{s}\omega})P(d\omega)ds\\
&=&\tilde{\mu}(A).
\end{eqnarray*}
Let $\pi_{\Omega}:\Omega\times E\rightarrow\Omega$, $\pi_{\Omega}(\omega,x)=\omega$, be the projection onto $\Omega$. By (\ref{4})  and the measure preserving property of $\{\theta_t\}$, we get $\tilde{\mu}\circ\pi^{-1}_{\Omega}=P$.
Hence $\tilde{\mu}$ is an invariant probability measure of $\Psi$ (cf. \cite[Definition 1.4.1]{A}).

By (\ref{zsw}) and the measure preserving property of $\{\theta_t\}$, we get
\begin{eqnarray*}
\tilde{\mu}(A)&=&\int_{\mathbb{R}}\int_{\Omega}\frac{\delta_s{([0,T\omega])}\cdot\delta_{Y(s,\omega)}(A_{\theta_{s}\omega})}{T\omega}P(d\omega)ds\\
&=&\int_{\mathbb{R}}\int_{\Omega}\frac{\delta_s{([0,T\omega])}\cdot\delta_{Y(s,\theta_{-s}\omega)}(A_{\omega})}{T\omega}P(d\omega)ds\\
&=&\int_{\Omega}\frac{1}{T\omega}\int_0^{T\omega}\delta_{Y(s,\theta_{-s}\omega)}(A_{\omega})dsP(d\omega).
\end{eqnarray*}
Then, the random factorization of $\tilde{\mu}$ is given by
$$
(\tilde{\mu})_{\omega}=\frac{1}{T\omega}\int_0^{T\omega}\delta_{{Y}(s,\theta_{-s}\omega)}ds,
$$
which is supported by $L^{\omega}$. Therefore, the proof is complete. \hfill\fbox\\

We now consider weak random periodic measures and invariant measures of a semi-flow $\varphi$. Define $\overline{E}:=\mathbb{R}\times E$ and
$$
\overline{\Psi}(t,\omega)(s,x)=(t+s,\varphi(t+s,s,\theta_{-s}\omega)x),\ \ \ \ \omega\in\Omega,\,s\in\mathbb{R},\, x\in E,\, t\ge0.
$$
Then, $\overline{\Psi}:[0,\infty)\times\Omega\times \overline{E}\rightarrow \overline{E}$ is a measurable random dynamical
system on $(\overline{E},{\cal E}(\overline{E}))$ over the metric dynamical system $(\Omega,{\cal F},P,(\theta_t)_{t\in\mathbb{R}})$. Assume that  $\varphi$ has a weak random periodic solution $(Y,T)$. Define
$$
\overline{Y}(s,\omega)=(s,Y(s,\omega)),\ \ \ \ \omega\in\Omega,\,s\in\mathbb{R},
$$
and
$$
\eta_t(s,x)=(t+s,x),\ \ \ \ \omega\in \Omega,\, s,t\in\mathbb{R},\, x\in E.
$$
Then, by (\ref{3333}), we obtain that for almost all $\omega\in\Omega$,
\begin{eqnarray}\label{July24}
   \overline{\Psi}(t,\theta_s\omega)\overline{Y}(s,\omega)=\overline{Y}(t+s,\omega),\ \ \overline{Y}(s+T\omega,\theta_{-T\omega}\omega)=\eta_{T\omega}\circ\overline{Y}(s,\omega),\ \ s\in\mathbb{R},\, t\ge0.\ \ \ \
\end{eqnarray}

Define
$$
\overline{\Upsilon}_t(\omega,s,x)=(\theta_t\omega,\overline{\Psi}(t,\omega)(s,x)),\ \ \ \ \omega\in\Omega,\,s\in\mathbb{R},\, x\in E,\, t\ge0,
$$
and
$$
\tilde{\eta}_t(\omega,s,x)=(\omega,t+s,x),\ \ \ \ \omega\in\Omega,\,s,t\in\mathbb{R},\, x\in E.
$$
Denote by ${\cal P}(\Omega\times \overline{E})$ the set of all probability measures on $(\Omega\times \overline{E}, \cal F\otimes\cal B(\overline{E}))$. Let $\pi_{\Omega}:\Omega\times \overline{E}\rightarrow\Omega$, $\pi_{\Omega}(\omega,(s,x))=\omega$, be the projection onto $\Omega$.

\begin{defi}
    A weak random periodic probability measure of $\varphi$ is a pair of measurable maps $\mu: \mathbb{R}\times\Omega\rightarrow{\cal P}(\Omega\times \overline{E})$ and $T:\Omega\rightarrow(0,\infty)$ such that for almost all $\omega\in\Omega$,
    $$
    \overline{\Upsilon}_t\mu(s,\omega)=\mu({t+s},\omega), \ \
    \mu({s+T\omega},\theta_{-T\omega}\omega)\circ \tilde{\eta}_{T\omega}=\mu(s,\omega), \ \ \ \ \forall s\in\mathbb{R},\, t\ge0.
    $$
\end{defi}

\begin{thm}\label{semi}
    If a stochastic semi-flow $\varphi:\Delta\times\Omega\times E \rightarrow E$ has a weak random periodic solution  $Y:\mathbb{R}\times\Omega\rightarrow E$
    and $T:\Omega\rightarrow(0,\infty)$, then it has a weak random periodic probability measure. If in addition (\ref{zsw}) holds for almost all $\omega\in\Omega$, then there exists a weak-invariant probability measure $\tilde{\mu}$ on $\cal F\otimes\cal B(\overline{E})$ satisfying $\tilde{\mu}\circ\pi^{-1}_{\Omega}=P$,
$$
\overline{\Upsilon}_t\tilde{\mu}(A)=\tilde{\mu}(A),\ \ \ \ \forall A\in \cal F\otimes\{\emptyset,\mathbb{R}\}\times \cal B({E}),\,t\ge 0,
$$
and its random factorization is supported by
$$
L^{\omega}:=\{\overline{Y}(s,\theta_{-s}\omega):s\in[0,T\omega)\}.
$$
\end{thm}
\noindent {\bf Proof.}\ \ For $s\in\mathbb{R}$ and $\omega\in\Omega$, define
\begin{eqnarray*}
\mu(s,\omega)(A)=\delta_{\overline{Y}(s,\omega)}(A_{\theta_{s}\omega}),\ \ \ \ A\in \cal F\otimes\cal B(\overline{E}),
\end{eqnarray*}
where $A_{\omega}$ is the $\omega$-section of $A$. Then,  $\mu(s,\omega)\in {\cal P}(\Omega\times \overline{E})$.

We have
\begin{eqnarray*}
(\overline{\Upsilon}^{-1}_t(A))_{\omega}&=&\{(s,x):(\theta_t\omega,\overline{\Psi}(t,\omega)(s,x))\in A\}\\
&=&\{(s,x):\overline{\Psi}(t,\omega)(s,x)\in A_{\theta_t\omega}\}\\
&=&\overline{\Psi}^{-1}(t,\omega)A_{\theta_t\omega}.
\end{eqnarray*}
Then, by (\ref{July24}), we obtain that for almost all $\omega\in\Omega$,
\begin{eqnarray*}
{\overline{\Upsilon}}_t\mu(s,\omega)(A)
&=&\mu(s,\omega)({\overline{\Upsilon}}^{-1}_t(A))\nonumber\\
&=&\delta_{\overline{Y}(s,\omega)}(({\Upsilon}^{-1}_t(A))_{\theta_{s}\omega})\nonumber\\
&=&\delta_{\overline{\Psi}(t,\theta_{s}\omega)\overline{Y}(s,\omega)}(A_{\theta_{t+s}\omega})\nonumber\\
&=&\delta_{\overline{Y}(t+s,\omega)}(A_{\theta_{t+s}\omega})\nonumber\\
&=&\mu({t+s},\omega)(A),
\end{eqnarray*}
and
\begin{eqnarray*}
\mu({s+T\omega},\theta_{-T\omega}\omega)(\tilde{\eta}_{T\omega} A)
&=&\delta_{\overline{Y}({s+T\omega},\theta_{-T\omega}\omega)}((\tilde{\eta}_{T\omega}A)_{\theta_{s}\omega})\nonumber\\
&=&\delta_{\eta_{T\omega}\circ \overline{Y}(s,\omega)}(\eta_{T\omega}A_{\theta_{s}\omega})\nonumber\\
&=&\delta_{\overline{Y}(s,\omega)}(A_{\theta_{s}\omega})\nonumber\\
&=&\mu({s},\omega)(A).
\end{eqnarray*}
Thus, $\mu$ is a weak random periodic probability measure of
$\varphi$.

For $A\in \cal F\otimes\cal B(\overline{E})$, define
$$
\tilde{\mu}(A):=\int_{\Omega}\frac{1}{T\omega}\int_0^{T\omega}\mu(s,\omega)(A)dsP(d\omega).
$$
Then, by using (\ref{July24}) and following the same argument of the proof of Theorem \ref{invar}, we can complete the proof.\hfill\fbox\\

\section{Examples}\setcounter{equation}{0}

In this section, we use examples  to illustrate the weak random
periodic phenomena of dynamical systems induced by random and stochastic differential equations.

First, we investigate the periodic behavior of RDEs of type (\ref{KKLL}) by virtue of weak random periodic solutions.

\begin{exa}\label{41}  Let $X(\varpi)$ be  a positive random variable and $a_k,b_k\in\mathbb{R}$, $1\le k\le N$, for some $N\in\mathbb{N}$.
Consider the following RDE:
\begin{equation}\label{KKLLL}
\frac{d^2x(t)}{dt^2}=\sum_{k=1}^N[a_k\sin(kX(\varpi)(t+s))+b_k\cos(kX(\varpi)(t+s))],
\end{equation}
where $s\in\mathbb{R}$. Note that Equation (\ref{KKLLL}) is equivalent to
\begin{eqnarray}\label{RDEE}
\left\{
\begin{array}{l}
dx_1(t)=x_2(t)dt,\\
dx_2(t)=\left\{\sum_{k=1}^N[a_k\sin(kX(\varpi)(t+s))+b_k\cos(kX(\varpi)(t+s))]\right\}dt.
\end{array}
\right.
\end{eqnarray}

Denote by $\nu$ the distribution of $X(\varpi)$. Define
$$
V=\left\{(x,y):x\in(0,\infty),\,y\in\left[0,\frac{2\pi}{x}\right)\right\}.
$$
We equip $(V,{\cal B}(V))$ with the probability measure $P_V$:
$$
P_V(A)=\int_0^{\infty}\int_0^{\infty}\frac{x1_A(x,y)}{2\pi}dy\nu(dx),\ \ \ \ A\in {\cal B}(V).
$$
Define
$$
g_{x,y}(t)=\sum_{k=1}^N[a_k\sin(kx(t+y))+b_k\cos(kx(t+y))],\ \ \ \ t\in \mathbb{R},\,(x,y)\in V,
$$
and
$$
\Omega=\{g_{x,y}:(x,y)\in V\}.
$$
Set $J:V\mapsto \Omega, J(x,y)=g_{x,y}$. Define
$$
{\cal F}=J({\cal B}(V)),\ \ P=P_V\circ J^{-1},
$$
and
$$
\theta_t\omega(s)=\omega(t+s),\ \ \ \ \omega\in\Omega,\,s,t\in\mathbb{R}.
$$
Then, $(\Omega,{\cal F},P,(\theta_t)_{t\in\mathbb{R}})$ is a metric dynamical system and Equation (\ref{RDEE}) is equivalent to the following RDE:
\begin{eqnarray}\label{RDEEP}
\left\{
\begin{array}{l}
dx_1(t)=x_2(t)dt,\\
dx_2(t)=\omega(t)dt.
\end{array}
\right.
\end{eqnarray}

The random dynamical system $\Psi:\mathbb{R}\times
\Omega\times \mathbb{R}^2\rightarrow \mathbb{R}^2$ induced by Equation (\ref{RDEEP}) is given by
\begin{eqnarray*}
&&\Psi(t,g_{x,y})(x_1,x_2)\\
&=&\left( x_1+x_2t-\sum_{k=1}^N\frac{a_k[\sin(kx (t+y) )-kxt\cos(kxy)-\sin(kxy)]}{k^2x^2}\right.\\
&&\ \ \ \ \ \ -\sum_{k=1}^N\frac{b_k[\cos(kx (t+y) )+kxt\sin(kxy)-\cos(kxy)]}{k^2x^2},\\
&&\left.\ \ \ x_2+\sum_{k=1}^N\frac{-a_k[\cos(kx (t+y) )-\cos(kxy)]+b_k[\sin(kx (t+y) )-\sin(kxy)]}{kx}\right),
\end{eqnarray*}
where $t\in \mathbb{R}$, $(x,y)\in V$ and $(x_1,x_2)\in\mathbb{R}^2$. Fix a ${\cal B}(\mathbb{R})$-measurable function $h:\mathbb{R}\rightarrow\mathbb{R}$. For $t\in \mathbb{R}$ and $(x,y)\in V$, define
\begin{eqnarray}\label{uni}
&&Y(t,g_{x,y})\nonumber\\
&=&\left(h(x) -\sum_{k=1}^N\frac{a_k\sin(kx (t+y) )+b_k\cos(kx (t+y) )}{k^2x^2},\right.\nonumber\\
&&\left.\ \ \ \sum_{k=1}^N\frac{-a_k\cos(kx (t+y) )+b_k\sin(kx (t+y) )}{kx}\right),
\end{eqnarray}
and
$$
Tg_{x,y}=\frac{2\pi}{x}.
$$
Then, $(Y,T)$ is a weak random periodic solution of $\Psi$. Further,  by Theorem \ref{invar}, we know that $\Psi$ has an invariant probability measure.
\end{exa}

\begin{rem} Theorem \ref{invar} shows that, if a random dynamical system has a weak random periodic solution, then it has an invariant probability measure induced by this solution. On the other hand, Example \ref{41} shows that different weak random periodic solutions can be obtained for the same random dynamical system through choosing different functions $h$ in (\ref{uni}). Therefore, it is interesting to consider how weak random periodic solutions affect the ergodicity of random dynamical systems.
\end{rem}

Next, we consider a system of  RDEs driven by periodic multiplicative noises.

\begin{exa}\label{exa1122} Suppose $d\ge 1$. Denote by $C(\mathbb{R};\mathbb{R}^d)$ and $C^1(\mathbb{R};\mathbb{R}^d)$ the spaces of all continuous and continuously differentiable $\mathbb{R}^d$-valued functions  on $\mathbb{R}$, respectively. We equip $C(\mathbb{R};\mathbb{R}^d)$ with topology of locally uniform convergence. Define $\Omega=C^1(\mathbb{R};\mathbb{R}^d)$ and
$$
{\cal F}=\{A\cap \Omega:A\in{\cal B}(C(\mathbb{R};\mathbb{R}^d))\}.
$$
 For $\omega=(\omega_1,\dots,\omega_d)\in \Omega$ and $s,t\in \mathbb{R}$, set $(\theta_t\omega)(s)= \omega(t + s)-\omega(t)$.

We choose $\omega^1,\omega^2,\dots\in \Omega$ with periods $T_1<T_2<\cdots$, respectively,  and $a_1,a_2,\dots\in(0,\infty)$ satisfying $\sum_{n=1}^{\infty}a_n=1$.  Define
  $$
 \Omega_n:=\{\theta_t\omega^n:0\le t<T_n\},\ \ \ \ n\in\mathbb{N}.
  $$
Denote by ${\cal L}$ the Lebesgue measure on $\mathbb{R}$. We define  a probability measure $P$ on $(\Omega,{\cal F})$ by
  $$
  P\left(\Omega\setminus \bigcup_{n=1}^{\infty}\Omega_n\right)=0,
  $$
  and
  $$
  P(\{\theta_t\omega_n:t\in A\})=\frac{a_n{\cal L}(A)}{T_n},\ \ \ \ \forall A\in {\cal B}([0,T_n)),\,n\in\mathbb{N}.
  $$
Set
  $$
  T\omega=T_n,\ \ \ \ \forall \omega\in \Omega_n,\,n\in\mathbb{N},
  $$
and
$$
T\omega=1,\ \ \ \ \forall\omega\notin\bigcup_{n=1}^{\infty}\Omega_n.
$$
  Then, $\{\theta_t\}$ are $P$-measure preserving and $T:\Omega\rightarrow(0,\infty)$ is a measurable map such that for almost all $\omega\in\Omega$,
$$
\omega(s+T\omega)=\omega(s),\ \ \ \ \forall s\in\mathbb{R}.
$$

Let $A$ be a $d\times d$ hyperbolic matrix and $\sigma=(\sigma_{ij})_{1\le i,j\le d}$ with $\sigma_{ij}:\mathbb{R}^d\rightarrow \mathbb{R}$ being Lipschitz-continuous and satisfying
$\sigma_{ij}(x)=o(|x|)$ as $|x|\rightarrow\infty$. Consider the RDE
    \begin{eqnarray}\label{zh11x}
        dx(t)=Ax(t)dt+\sigma(x(t))d\omega(t).
    \end{eqnarray}
By \cite[Theorem 22.1]{A1}, we know that  (\ref{zh11x}) has a $T_n$-periodic solution  for each $\omega^n$, which is denoted by $x(t,\omega^n)$. For $s,t\in\mathbb{R}$, define
$$
Y(t,\theta_s\omega^n)=x(t+s,\omega^n).
$$
Then, $(Y,T)$ is a weak random periodic solution of the random dynamical system induced by Equation (\ref{zh11x}). Further, by Theorem \ref{invar}, we know that  this random dynamical system  has an invariant probability measure.
\end{exa}

The third example is concerned with a random dynamical system induced by stochastic differential equations (SDEs), which is is an extension of the example given by Zhao and Zheng (see \cite[Section 2]{ZZ2}).

\begin{exa}\label{EK} Let $\Omega:=C(\mathbb{R};\mathbb{R})$ and $\{\omega(t)\}_{t\in\mathbb{R}}$ be a one-dimensional two-sided Brownian motion on the path space $(\Omega, {\cal B}(\Omega),P)$
with  $\theta$ being the shift operator $(\theta_t\omega)(s)= \omega(t + s)-\omega(t)$ for $s,t\in \mathbb{R}$. We define an equivalence relation $\sim$ on
$\Omega$ by $\omega\sim \omega'$ if and only
if there exists $t\in\mathbb{R}$ such that
$\omega'=\theta_t\omega$. Denote by $\Omega':=\Omega/\sim$ the
quotient space of $\Omega$.

Before stating the example, we present a proposition. To the best of our knowledge, this is a novel result in the literature, which is of independent interest.

  \begin{pro}\label{wiener}
  $(\Omega',{\cal B}(\Omega'),P')$ is isomorphic (mod  0) to $[0,1]$ with the Lebesgue measure.
  \end{pro}

\noindent {\bf Proof.}\ \   Since $\Omega$ is  countably generated, $\Omega'$ is also countably generated.  Then,  $(\Omega',{\cal B}(\Omega'),P')$ is a standard probability space (cf. \cite{M,R}).
  To show that $(\Omega',{\cal B}(\Omega'),P')$  is isomorphic (mod  0) to $[0,1]$ with the Lebesgue measure, it is sufficient to prove that the probability
  space $(\Omega',{\cal B}(\Omega'),P')$ has no atom, equivalently, $P([w])=0$ for almost all
  $\omega\in\Omega$.  Hereafter, $[\omega]$ denotes the equivalence class of $\omega$. Note that
  $$[w]=\{\theta_t\omega:t\in[0,\infty)\}\cup\{\theta_t\omega:t\in(-\infty,0]\}.
  $$ By symmetry,  it is sufficient to  show that for almost all $\omega\in\Omega$,
  $$P^*(\{\theta_t\omega|_{[0,\infty)}:t\in[0,\infty)\})=0,
  $$
   where $P^*$ is the restriction of $P$ on $C[0,\infty)$.

Define
$$
A=\left\{\omega\in C[0,\infty):
\limsup_{t\rightarrow\infty}\omega(t)=\infty\ {\rm and}\
\liminf_{t\rightarrow\infty}\omega(t)=-\infty\right\}. $$ It is
well-known that $P^*(A)=1$. Fix an $\omega_0\in A$. We choose
$0<t_1<t_2<\infty$ such that $\omega_0(t_1)=\omega_0(t_2)=0$ and
$$
m_1:=\max_{0\le u\le t_1}\omega_0(u)>0,\ \ m_2:=\max_{t_1\le u\le
t_2}\omega_0(u)>0.
$$
Then, there exists $\varepsilon>0$ such that
$$
m_1=\max_{\varepsilon\le u\le t_1+\varepsilon}\omega_0(u),\ \
m_2=\max_{t_1+\varepsilon\le u\le t_2+\varepsilon}\omega_0(u).
$$
Define $$B=\left\{\omega\in C[0,\infty): \left[\max_{t_1\le u\le
t_2}\omega(u)\right]-\left[\max_{0\le u\le
t_1}\omega(u)\right]=m_2-m_1\right\}.
$$
Then,
\begin{equation}\label{local}
\{\theta_t\omega_0:t\in[0,\varepsilon)\}\subset B.
\end{equation}

For $\omega\in C[0,\infty)$, we have
\begin{eqnarray*}
&&\left[\max_{t_1\le u\le t_2}\omega(u)\right]-\left[\max_{0\le u\le
t_1}\omega(u)\right]\\
&=&\left[\max_{t_1\le u\le
t_2}\{\omega(u)-\omega(t_1)\}\right]-\left[\max_{0\le u\le
t_1}\{\omega(u)-\omega(t_1)\}\right]\\
&=&\left[\max_{0\le v\le
t_2-t_1}\{\omega(t_1+v)-\omega(t_1)\}\right]-\left[\max_{0\le v\le
t_1}\{\omega(t_1-v)-\omega(t_1)\}\right]\\
&:=&M_2-M_1.
\end{eqnarray*}
It is known that $M_2$ and $M_1$ are two independent continuous
random variables such that (cf. \cite[page 96]{KS})
$$
P^*(M_2\in dx)=\frac{2}{\sqrt{2\pi
(t_2-t_1)}}e^{-\frac{x^2}{2(t_2-t_1)}}dx,\ \ P^*(M_1\in
dx)=\frac{2}{\sqrt{2\pi t_1}}e^{-\frac{x^2}{2t_1}}dx;\ \ \ \ x>0.
$$
Then, $P^*(B)=0$, which together with (\ref{local}) implies that
$$P^*(\{\theta_t\omega_0:t\in[0,\varepsilon)\})=0.$$

Applying the similar argument, we can show that for any $s\ge 0$,
there exists $\varepsilon_{s}>0$ such that
$$P^*(\{\theta_t(\theta_s\omega_0):t\in[0,\varepsilon_s)\})=0,$$
which implies that
\begin{equation}\label{local2}
P^*(\{\theta_{t}\omega_0:t\in[s,s+\varepsilon_s)\})=0,\ \ \ \ \forall
s\ge0.
\end{equation}
Define
$$
C=\sup\{c: P^*(\{\theta_t\omega_0: t\in[0,c)\})=0\}.
$$
Then, by (\ref{local2}), we obtain that $C=\infty$. Hence
$$
P^*(\{\theta_t\omega_0: t\in[0, \infty)\})=0.
$$
Since $\omega_0\in A$ is arbitrary, we conclude that for almost all
$\omega\in\Omega$,
$$P^*(\{\theta_t\omega|_{[0,\infty)}:t\in[0,\infty)\})=0.
$$
Therefore, $(\Omega',{\cal B}(\Omega'),P')$ is isomorphic (mod  0)
to $[0,1]$ with the Lebesgue measure.\hfill\fbox\\

We consider the SDE
\begin{eqnarray}\label{zh11}
\left\{
\begin{array}{l}
dx(t)=\{x(t)-y(t)-x(t)[x^2(t)+y^2(t)]\}dt+x(t)\circ d\omega(t),\\
dy(t)=\{x(t)+y(t)-y(t)[x^2(t)+y^2(t)]\}dt+y(t)\circ d\omega(t),
\end{array}
\right.
\end{eqnarray}
where $\circ d\omega(t)$ denotes the Stratonovich stochastic integral. Using polar
coordinates
$$
x=\rho\cos(2\pi\alpha),\ \ y=\rho\sin(2\pi\alpha),
$$
we can transform Equation (\ref{zh11}) on $\mathbb{R}^2$ to the following equation on $[0,\infty)\times[0,\infty)$:
\begin{eqnarray}\label{zh22}
\left\{
\begin{array}{l}
d\rho(t)=[\rho(t)-\rho^3(t)]dt+\rho(t)\circ d\omega(t),\\
d\alpha(t)=\frac{1}{2\pi}dt.
\end{array}
\right.
\end{eqnarray}
Equation (\ref{zh22}) has a unique closed form solution as follows:
$$
\rho(t,\alpha_0,\rho_0,\omega)=\frac{\rho_0e^{t+\omega(t)}}{(1+2\rho_0^2\int_0^te^{2s+2\omega(s)}ds)^{\frac{1}{2}}},\ \ \alpha(t,\alpha_0,\rho_0,\omega)=\alpha_0+\frac{t}{2\pi}.
$$
We can check that
$$
\rho^*(\omega)=\left(2\int_{-\infty}^0e^{2s+2\omega(s)}ds\right)^{-\frac{1}{2}}
$$
is the stationary solution of the first equation of (\ref{zh22}), i.e.,
$$
\rho(t,\alpha_0,\rho^*(\omega),\omega)=\rho^*(\theta_t\omega).
$$

By Proposition \ref{wiener}, we can choose a surjective measurable map
$T':\Omega'\rightarrow (0,\infty)$.  Define $T\omega:=T'[\omega]$ for
$\omega\in \Omega$. Then, $T$ is a measurable map on $\Omega$. Define
$$
\Psi^*(t,\omega)(\alpha_0,\rho_0)=\left(\left(\alpha_0+\frac{t}{T\omega}\right)\ {\rm mod}\ 1,\ \rho(t,\alpha_0,\rho_0,\omega)\right).
$$
We find that
$$
\Psi^*(0,\omega)(\alpha,\rho)=(\alpha,\rho),\ \Psi^*(t+s,\omega)=\Psi^*(t,\theta_s\omega)\Psi^*(s,\omega),\ \ \ \ \forall (\alpha,\rho)\in [0,1)\times [0,\infty),\, s,t\in\mathbb{R}.
$$
Hence $\Psi^*(t,\omega)=(\Psi^*_1(t,\omega), \Psi^*_2(t,\omega))$ defines a random dynamical system on the cylinder $[0,1)\times [0,\infty)$.
Next we transform the random dynamical system $\Psi^*$
back to $\mathbb{R}^2$. For $(x,y)\in\mathbb{R}^2, x =\rho\cos(2\pi\alpha), y =\rho\sin(2\pi\alpha)$, define
$$
\Psi(t,\omega)(x,y)=\left(\Psi^*_2(t,\omega)(\alpha,\rho)\cdot\cos[2\pi\Psi^*_1(t,\omega)(\alpha,\rho)],\,
\Psi^*_2(t,\omega)(\alpha,\rho)\cdot\sin[2\pi\Psi^*_1(t,\omega)(\alpha,\rho)]\right).
$$

Now we investigate weak random periodic solutions of the random dynamical system $\Psi$. Fix an $\alpha_0\in [0,1)$ and define
$$
Y(t,\omega)=\left(\rho^*(\theta_t\omega)\cos\left[2\pi\alpha_0+\frac{2\pi t}{T\omega}\right],\,\rho^*(\theta_t\omega)\sin\left[2\pi\alpha_0+\frac{2\pi t}{T\omega}\right]\right).
$$
Then, we have
\begin{eqnarray*}
&&\Psi(T\omega,\theta_{-T\omega}\omega)Y_0(\theta_{-T\omega}\omega)\\
&=&\Psi(T\omega,\theta_{-T\omega}\omega)(\rho^*(\theta_{-T\omega}\omega)\cos(2\pi\alpha_0),
\rho^*(\theta_{-T\omega}\omega)\sin(2\pi\alpha_0))\\
&=&(\Psi_2^*(T\omega,\theta_{-T\omega}\omega)\cdot\cos[2\pi\Psi^*_1(T\omega,\theta_{-T\omega}\omega)(\alpha_0,\rho^*(\theta_{-T\omega}w))],\\
&&\ \,
\Psi_2^*(T\omega,\theta_{-T\omega}\omega)\cdot\sin[2\pi\Psi^*_1(T\omega,\theta_{-T\omega}\omega)(\alpha_0,\rho^*(\theta_{-T\omega}w))])\\
&=&\left(\rho(T\omega,\alpha_0,\rho^*(\theta_{-T\omega}w),
\theta_{-T\omega}\omega)\cdot\cos\left[2\pi\left(\left(\alpha_0+\frac{T\omega}{T(\theta_{-T\omega}\omega)}\right)\ \ {\rm
mod}\
1\right)\right],\right.\\
&&\ \ \left.\rho(T\omega,\alpha_0,\rho^*(\theta_{-T\omega}w),
\theta_{-T\omega}\omega)\cdot\sin\left[2\pi\left(\left(\alpha_0+\frac{T\omega}{T(\theta_{-T\omega}\omega)}\right)\ \ {\rm
mod}\
1\right)\right]\right)\\
&=&(\rho^*(\omega)\cos(2\pi\alpha_0),\,\rho^*(\omega)\sin(2\pi\alpha_0))\\
&=&Y_0(\omega),
\end{eqnarray*}
which implies that (\ref{26}) holds. Therefore, by
Proposition \ref{thm2}, we find that $(Y,T)$ is a weak
random periodic solution of $\Psi$. Further, by Theorem \ref{invar}, we conclude that $\Psi$ has an invariant probability measure.

\end{exa}

\begin{rem}
By Proposition \ref{wiener}, we know that there are different choices of the random period map $T'$. Hence, Example \ref{EK} implies that the weak random periodic solution of a
  random dynamical system is not necessary to be a random periodic solution defined as in \cite{FZ,ZZ2}.
  \end{rem}

The last example is related to SDEs in random environments.

\begin{exa}\label{EK} Let $Z$ be an $\mathbb{N}$-valued random variable on a probability space $(\Omega_1,{\cal F}_1,P_1)$. Suppose $d\ge 1$, $\Omega_2:=C(\mathbb{R};\mathbb{R}^d)$, and $\{\omega_2(t)\}_{t\in\mathbb{R}}$ be a $d$-dimensional two-sided Brownian motion on the path space $(\Omega_2, {\cal B}(\Omega_2),P_2)$
with  $\theta_2$ being the shift operator $(\theta_{2,t}\omega_2)(s)= \omega_2(t + s)-\omega_2(t)$ for $s,t\in \mathbb{R}$.
Let $m\in\mathbb{N}$, $T_n>0$, $b_{n,i}: \mathbb{R}\times \mathbb{R}^m\rightarrow \mathbb{R}$, $\sigma_{n,ij}: \mathbb{R}\times \mathbb{R}^m\rightarrow \mathbb{R}$, $1\le i\le m,\, 1\le j\le d$, $n\in\mathbb{N}$. Suppose
$$
b_{n,i}(t+T_n,x)=b_{n,i}(t,x),\ \ \sigma_{n,ij}(t+T_n,x)=\sigma_{n,ij}(t,x),\ \ \ \ t\in \mathbb{R},\, x\in \mathbb{R}^m,\,n\in\mathbb{N}.
$$
Assume that, for each $n\in\mathbb{N}$, the stochastic semi-flow $\varphi_n$ induced by the following SDE has a random periodic solution $\{Y_n(t,\omega_2)\}$ with period $T_n$:
$$
dx_n(t)=b_n(t,x_n(t))dt+\sigma_n(t,x_n(t))dw_2(t).
$$
We refer the reader to \cite{FZZ, FWZ3,FZ,D1, D2} for various concrete examples of random periodic solutions of
non-autonomous SDEs.

 Let $(\Omega,{\cal F},P)=(\Omega_1,{\cal F}_1,P_1)\times (\Omega_2, {\cal B}(\Omega_2),P_2)$ and $\theta_t(\omega_1,\omega_2)=(\omega_1,\theta_{2,t}(\omega_2))$. We consider the following SDE in random environment:
$$
dx(t)=b_Z(t,x(t))dt+\sigma_Z(t,x(t))dw_2(t),
$$
which induces a stochastic semi-flow $\varphi$:
$$
\varphi(t,s,\omega_1,\omega_2)=\varphi_{Z(\omega_1)}(t,s,\omega_2),\ \ \ \ (\omega_1,\omega_2)\in\Omega,\,s\le t.
$$
Define
$$
T(\omega_1,\omega_2)=T_{Z(\omega_1)},\ \ Y(t,\omega_1,\omega_2)=Y_{Z(\omega_1)}(t,\omega_2),\ \ \ \ (\omega_1,\omega_2)\in\Omega,\, t\in\mathbb{R}.
$$
Then, $(Y,T)$ is a weak random periodic solution of  $\varphi$, which has a weak-invariant probability measure by Theorem \ref{semi}.

\end{exa}


\vskip 1cm

{ \noindent {\bf\large Acknowledgements}\ \
W. Sun acknowledges the financial support of the Natural Sciences and Engineering Research Council of Canada. Z.H. Zheng acknowledges financial supports of the NSF of China (No. 12090014, 12031020, 11671382), CAS Key Project of Frontier Sciences (No. QYZDJ-SSW-JSC003), the Key Lab of Random Complex Structures and Data Sciences CAS and National Center for Mathematics and Interdisciplinary Sciences CAS.

\end{document}